\documentclass[9pt, a4paper]{amsart}
\usepackage{amssymb,amsmath,mathrsfs,amsthm,amscd}
\usepackage{epsf}

\newtheorem{theorem}{Theorem}
\newtheorem*{conjecture}{Toral rank conjecture, \cite{Halp}}

\newtheorem{lemma}[theorem]{Lemma}

\theoremstyle{definition}
\newtheorem{definition}[theorem]{Definition}

\newcommand{\R}{\mathbb{R}}
\newcommand{\Z}{\mathbb{Z}}
\newcommand{\Cbb}{\mathbb{C}}

\renewcommand{\ge}{\geqslant}
\newcommand{\Ps}{\mathcal{P}}
\renewcommand{\a}{\alpha}
\newcommand{\pd}{\partial}
\newcommand{\dt}{\pd_t}
\newcommand{\da}{\pd_{\a}}
\newcommand{\codim}{\mathop{\rm codim}}
\newcommand{\es}{\varnothing}

\begin{document}
\title{Doubling operation for polytopes and torus actions}
\author{Yury Ustinovsky}
\date{}
\maketitle

In this note we give the definition of the "doubling operation" for
simple polytopes, find the formula for the $h$-polynomial of new
polytope. %in terms of the initial.
As an application of this
operation we establish the relationship between moment-angle
manifolds and their real analogues and prove the toral rank
conjecture for moment-angle manifolds $Z_P$.

Let $P$ be a simple n-dimensional polytope with $m$ facets:
$$P=\{x\in \R^n\mid (x,a_i)+b_i\ge0, \quad a_i \in
\R^n, \ i=1\ldots{m}\}.$$ Then $P$ can be identified with the
intersection of the positive orthant $\R^m_{\ge}$ and the image of
the affine map $i_P: \R^n \to \R^m,$\quad $i_P(x)=A_P(x)+b_P$, where
$A_P$ is  $m\times n$ matrix with rows $a_i$ and $b_P$ is column
vector $(b_1,\dots, b_m)$. Suppose the image of the map $i_P$ is
specified by the system of $m-n$ equations in $\R^m$:
\begin{equation}
\sum\limits_{i=1}^m{c_{j,i}x_i=q_i},\ \ j=1,\dots,m-n.
\end{equation}
Then the image of the polytope $P$ under map $i_P$ is the set of
points $p\in\R^m_{\ge}$ satisfying this equation.
\begin{definition}
\emph{The double} of the polytope $P$ is the polytope $L(P)$, given
by the system of equations in $\R^{2m}_{\ge}:$
$$\sum_{i=1}^m{c_{j,i}(x_i+x_i')=q_i},\ \ j=1,\dots,m-n.$$
\end{definition}

It is easy to see that $L(P)$ is $(m+n)$-dimensional simple polytope
with $2m$ facets, and $L(P\times Q)=L(P)\times L(Q)$. If
$\{v_1,\dots,v_m\}$ is the Gale diagram of $P$ (see \cite{Z}), then
$\{v_1, v_1,\dots,v_m, v_m\}$ is the Gale diagram of $L(P)$.

Let $f_{n-i,i}$ denote the number of faces of dimension $n-i$. Then
define the $h$-polynomial by the equation $h(P)(\a,t)=f(P)(\a-t,t)$,
where $f(P)(\a,t)=\a^n+f_{n-1,1}\a^{n-1}t+\dots+f_{0,n}t^n.$

\begin{lemma}
$h(L(P))(\a,t)=\sum{(-1)^{\codim G}(\a t)^{\codim G}(\a+t)^{m-\codim
G}h(G)(\a,t)},$ where $G$ ranges over the set of the faces of the
polytope, including $P$ itself.
\end{lemma}

Let $(\Ps, d)$ denote the differential ring of simple polytopes,
\cite{Buh}. Multiplication in $\Ps$ is given by the product of
polytopes, the sum is formal, and the differential is given by the
formula $d(P)=\sum{F_i}$, where the right side is the sum of the
facetes of $P$. Let $h: \Ps\to\mathbb{Z}[\a,t]$ be the ring
homomorphism determined by the $h$-polynomial, then
$h(d(P))(\a,t)=(\da+\dt)h(P)(\a,t),$ \cite{Buh}. This fact allows to
reformulate lemma 2:
$$h(L(P))(\a,t)=(\a+t)^m\Bigl(\sum\limits_{i=0}^{\infty}{\cfrac{(-1)^i}{i!}(\a t)^i(\a+t)^{-i}(\da+\dt)^i}\Bigr)h(P)(\a,t).$$
The differential operator in brackets is multiplicative, so it is
determined by its value on generators $\a$ and $t$.
\begin{theorem}
$h(L(P))(\a,t)=(\alpha+t)^{m-n}h(P)(\a^2,t^2)$.
\end{theorem}

Now let us consider the example of the application of the "doubling
operation".
\begin{definition}
Let $\mu: \Cbb^m\to\R^m_{\ge},\ \
\mu(z_1,\dots,z_m)=(|z_1|^2,\dots,|z_m|^2)$ be the projection on the
orbit space of the standard action of $T^m$. \emph{Moment-angle
manifold} $Z_P$ is the set $\mu^{-1}(i_P(P))$, (see. \cite{BPS}).
\end{definition}
\begin{definition}
Let $\pi: \R^m\to\R^m_{\ge},\ \
\pi(x_1,\dots,x_m)=(x_1^2,\dots,x_m^2)$ be the projection on the
orbit space of the standard action of $Z_2^m$. \emph{Real
moment-angle manifold} $R_P$ is the set $\pi^{-1}(i_P(P))$.
\end{definition}

There is natural action $T^m\colon Z_P$ and the orbit space of this
action can be identified with $P$. Similarly there is action
$\Z_2^m\colon R_P$. The simpleness of $P$ implies the existence of
the structure of a smooth manifold on the spaces $Z_P$ and $R_P$,
\cite{BPS}.

\begin{conjecture}
Let $X$ be a finite-dimensional topological space. If $T^r$ acts on
$X$ almost freely, then $\emph{hrk}(X)\ge 2^r$, where
$\emph{hrk}(X)=\sum_i{\dim{H^i(X)}}$.
\end{conjecture}

Manifolds $Z_P$ provide big class of spaces with the action of
$m$-dimensional torus, moreover the rank of subtorus in $T^m$ which
acts almost freely equals $m-s$ (see \cite{BPF}). We show that
$\mbox{hrk}(Z_P)\ge 2^{m-n}$, thus the toral rank conjecture holds
for the moment-angle manifolds.

\begin{lemma}\label{RMM}
Any moment-angle manifold $Z_P$ is homeomorphic to the real
moment-angle manifold $R_P$.
\end{lemma}

\begin{proof}
Assume $P$ is given by a system of the form (1). Then the
definitions imply that the manifold $Z_P$ is given by the system of
the equations in $\Cbb^m$: $$\sum_{i=1}^m{c_{j,i}|z_i|^2=q_i},\ \
j=1,\dots,m-n.$$ In the real coordinates $z_i=x_i+\sqrt{-1}y_i$ we
obtain the system:
$$\sum_{i=1}^m{c_{j,i}(x_i^2+y_i^2)=q_i},\
\ j=1,\dots,m-n, $$ which defines the manifold $R_{L(P)}$ in
$\R^{2m}$.
\end{proof}

\begin{theorem}\label{main}
The cohomology of an arbitrary real moment-angle manifold satisfy
inequality $\emph{hrk}(R_P)\ge 2^{m-n}$.
\end{theorem}
The next lemma is a consequence of the Mayer Vietoris long exact
sequence.
\begin{lemma}\label{HomolRk}
Let $M$ be compact manifold with boundary. Denote by $X$ the
manifold obtained from the two copies of $M$ by gluing them together
along common boundary. Then $\emph{hrk}(X)\ge \emph{hrk}(\pd M)$.
\end{lemma}

\begin{proof}[Proof of the theorem 7.] The case $n=1$ is trivial.
Suppose the statement of the theorem is proved for the polytopes of
the dimension $n-1$. Let $P$ be an arbitrary $n$-dimensional
polytope. Let $\R^m_+=\{\vec{x}\in \R^m|\ x_1\ge 0\}$ and denote by
$M$ the manifold $M=R_P\cap \R^m_+$. Then real moment-angle manifold
is obtained from the two copies of $M$: $R_P=M\cup_{\pd M}M$. This
decomposition satisfies lemma \ref{HomolRk}, thus
$\mbox{hrk}(R_P)\ge \mbox{hrk}(\pd M).$

Let us describe $\pd M$. By the definition of $M$, its boundary is
$R_P\cap\{\vec{x}\in \R^m|\ x_1= 0\}=\pi^{-1}(i_P(F_1))$~---
preimage of the facet $F_1$, given by the vanishing of the first
coordinate in $\R^m_{\ge}$, under the canonical projection on the
orbit space. Hence $\pi^{-1}(i_P(F_1))=R_{F_1}\times\Z_2^k$, where
$k$ is the number facets $F$ of $P$ such that $F_1\cap F=\es$ (thus,
the number of facets of $F_1$ equals $m-k-1$). By the assumption of
induction $\mbox{hrk}({\pd M})=2^k\mbox{hrk}(R_{F_1})\ge
2^k\cdot2^{(m-k-1)-(n-1)}.$ Finally we have: $\mbox{hrk}(R_P)\ge
\mbox{hrk}(\pd M)\ge 2^{m-n}.$
\end{proof}

Applying proved theorem and lemma \ref{RMM} we obtain the following
estimation for the cohomology rank of the moment-angle manifolds:
$\mbox{hrk}(Z_P)=\mbox{hrk}(R_{L(P)})\ge2^{2m-(m+n)}=2^{m-n},$ which
implies the statement of the toral rank conjecture for the
moment-angle manifolds.

It should be noted, that there is a generalization of the above
constructions for the simplicial complexes, which allows to prove
the toral rank conjecture for the class of moment-angle complexes,
see definitions in \cite{BPF}.

The author is grateful to T.\,E.~Panov and S.~Gitler for encouraging
to consider the "doubling operation". I would also like to thank
V.\,M.~Buchstaber and scientific adviser T.\,E.~Panov for suggesting
the problem and attention to the research.

\end{document}